\documentclass[a4paper]{amsart}  

\usepackage{amssymb} \usepackage{amscd} \usepackage{times}

\def\zbb{\mathbb{Z}}  
  
  \def\phi{\varphi}
 \def\p1{{\mathbb{P}^1_\zbb}}

\begin{document}

\title{ An Uniforme Estimate for Scalar Curvature Equation on Manifolds of dimension 4.}

\author{Samy Skander Bahoura}

\address{6, rue Ferdinand Flocon, 75018 Paris, France.}
              
\email{samybahoura@yahoo.fr, bahoura@ccr.jussieu.fr} 

\date{}

\maketitle

\begin{abstract}

We give an a priori estimate for the solutions of the prescribed scalar curvature equation on manifolds of dimension 4. We have an idea on the supremum of the solutions if we control their infimum.

\end{abstract}

\bigskip

\begin{center}

1. INTRODUCTION AND RESULT.

\end{center}

\bigskip

In this paper, we are on Riemannian manifold of dimension 4, $ (M,g) $ (not necessarily compact). Here we denote by $ \Delta_g=-\nabla^i(\nabla_i) $ the geometric Laplacian.

\bigskip

Let us consider the prescribed scalar curvature equation in four dimension:

$$ \Delta_g u+ R_g u =V u^3, \,\,\, u>0 \qquad (E) $$

where $ R_g $ is a scalar curvature of $ (M,g) $ and $ V $ the prescribed scalar curvature.

\bigskip

We assume:

$$ 0 < a \leq V(x) \leq b \,\,\, {\rm and} \,\,\, ||\nabla V||_{L^{\infty}(M)} \leq A \qquad (C). $$

In this paper, we want to prove an uniform estimate for the solutions of the eqution $ (E) $ with minimal conditions on the prescribed scalar curvature equation. Conditions like $ (C) $ are minimal.

\bigskip

Note that the equation $ (E) $ was studied when $ M = \Omega $ is a open set of $ {\mathbb R}^4 $, see for example, [B], [C-L] and when $ \Omega = {\mathbb S}_4 $ the unit sphere of dimension 4 by Li, see [L].

\bigskip

If we suppose $ V \in C^2(\Omega) $, Chen and Lin gave a $ \sup \times \inf $ inequality for the solutions of the equation $ (E) $. In [L], on the fourth unit sphere, Li study the same equation with the same conditions on $ V $, he obtains the boundedness of the energy and an upper bound for the product $ \sup \times \inf $. He use the simple blow-up analysis (for the definition of simple blow up points see for exemple [L]).

\bigskip
 
In [B], we can see (on a bounded domain of $ {\mathbb R}^4 $) that we have an uniform estimate for the solutions of the equation $ (E) $ if we control the infimum of those functions, with only Lipschitzian assumption on the prescribed scalar curvature $ V $.

\bigskip

Here we extend the result of [B], to general manifolds of dimension 4.

\bigskip

Note, if we assume $ V \equiv 1 $, Li and Zhang (see [L-Z 1]), have proved a $ \sup \times \inf $ inequality for the solutions of $ (E) $ on any Riemannian manifold of dimension 4.

\bigskip

If we suppose $ M $ compact, the existence result for this equation when $ V \equiv 1 $ was proved by T. Aubin (non conformally flat case and $ n\geq 6 $ ) and R. Schoen ( conformally flat case and $ n= 3, 4, 5 $). The previous equation with $ V \equiv 1 $ is called {\it the Yamabe equation}.

\bigskip

Note that, in diemsions $ n=3 $ and $ n\geq 5 $, we have many results about prescribed scalar curvature equation, see for example [B], [C-L], [L], [L-Z 1] and [L-Zh].

\bigskip

For example ( when $ M $ is compact), in [L-Zh], Li and Zhu have proved the compactness of the solutions of the Yamabe equation with the positive mass theorem. They also describe the blow-up points of the solutions ( only simple blow-up points). In [D], [L-Z 2] and [M], Druet, Li, Zhang and Marques have obtained the same result for the dimensions 4, 5, 6 and 7.

\bigskip

About the compactness of the solutions of the Yamabe equation, we can find in [L-Z 2] some conditions on the Weyl tensor to have this result. In [Au 2], T. Aubin have proved recently, the compactness of the soltuions of the Yamabe problem without other assumptions.

\bigskip

Note that here we have no assumption on energy. There is many results if we suppose the energy bounded. In our work, we use, in particular, the moving-plane method. This strong method was developped by Gidas-Ni-Nirenberg, see [G-N-N]. This method was used by many author to obtain uniform estimates, in dimension 2, see for example [B-L-S], in diemsnion greater than 3, see for example, [B], [ C-L], [L-Z 1] and [ L-Z 2].

\smallskip

We have:

\bigskip

{\it {\bf Theorem.} For all $ a, b, m >0 $, $ A \geq 0 $ with $ A \to 0 $ and all compact $ K $ of $ M $, there is a positive constant $ c=c(a, b, m, A, K, M, g) $ such that:

$$ \sup_K u \leq c \,\,\, {\rm if} \,\,\, \inf_M u \geq m, $$

for all solution $ u $ of $ (E) $ relatively to $ V $ with the conditions $ (C) $.} 

\bigskip

\bigskip

\bigskip

\begin{center}

2. PROOF OF THE THEOREM.

\end{center}

\bigskip

Let $ x_0 $ be a point of $ M $. We want to prove an uniform estimate around $ x_0 $.

\bigskip

Let $ (u_i)_i $ be a sequence of solutions of:

$$ \Delta u_i + R_g u_i =V_i{u_i}^3, \,\, u_i>0, $$

where $ V_i $ is such that:

$$ 0 < a \leq V_i(x) \leq b \,\,\, {\rm and} \,\,\, ||\nabla V_i||_{L^{\infty}(M)} \leq A_i \,\,\, {\rm with }\,\,\, A_i \to 0. $$

We argue by contradiction, we assume that the $ \sup $ is not bounded.

\bigskip

$ \forall \,\, c,R >0 \,\, \exists \,\, u_{c,R} $ solution de $ (E) $ telle que:

$$ R^2 \sup_{B(x_0,R)} u_{c,R} \geq c, \qquad (H) $$

\underbar {\bf Proposition 1:}{\it (blow-up analysis)} 

\smallskip

There is a sequence of points $ (y_i)_i $, $ y_i \to x_0 $  and two sequences of positive real numbers $ (l_i)_i, (L_i)_i $, $ l_i \to 0 $, $ L_i \to +\infty $, such that if we set $ v_i(y)=\dfrac{u_i[\exp_{y_i}(y/[u_i(y_i)])]}{u_i(y_i)} $, we have:

$$ 0 < v_i(y) \leq  \beta_i \leq 2, \,\, \beta_i \to 1. $$

$$  v_i(y)  \to \dfrac{1}{1+{|y|^2}}, \,\, {\rm uniformly \,\, on \,\, compact \,\, sets \,\, of } \,\, {\mathbb R}^4 . $$

$$ l_i u_i(y_i) \to +\infty. $$

\underbar {\bf Proof of the proposition 1:}

\bigskip

We use the hypothesis $ (H) $, we take two sequences $ R_i>0, R_i \to 0 $ and $ c_i \to +\infty $, such that,

$$ {R_i}^2 \sup_{B(x_0,R_i)} u_i\geq c_i \to +\infty, $$

Let, $ x_i \in  { B(x_0,R_i)} $, such that $ \sup_{B(x_0,R_i)} u_i=u_i(x_i) $ and $ s_i(x)=[R_i-d(x,x_i)] u_i(x), x\in B(x_i, R_i) $. Then, $ x_i \to x_0 $.

\bigskip

We have:

$$ \max_{B(x_i,R_i)} s_i(x)=s_i(y_i) \geq s_i(x_i)={R_i} u_i(x_i)\geq \sqrt {c_i}  \to + \infty. $$ 

We set :

$$ l_i=R_i-d(y_i,x_i),\,\, \bar u_i(y)= u_i [\exp_{y_i}(y)],\,\,  v_i(z)=\dfrac{u_i [ \exp_{y_i}\left ( z/[u_i(y_i)] \right )] } {u_i(y_i)}. $$

Clearly, we have, $ y_i \to x_0 $. We obtain:

$$ L_i= \dfrac{l_i}{(c_i)^{1/4}} [u_i(y_i)]=\dfrac{[s_i(y_i)]}{c_i^{1/4}}\geq \dfrac{c_i^{1/2}}{c_i^{1/4}}=c_i^{1/4}\to +\infty. $$

\bigskip

If $ |z|\leq L_i $, then $ y=\exp_{y_i}[z/ [u_i(y_i)]] \in B(y_i,\delta_i l_i) $ with $ \delta_i=\dfrac{1}{(c_i)^{1/4}} $ and $ d(y,y_i) < R_i-d(y_i,x_i) $, thus, $ d(y, x_i) < R_i $ and, $ s_i(y)\leq s_i(y_i) $. We can write,

$$ u_i(y) [R_i-d(y,y_i)] \leq u_i(y_i) l_i. $$

But, $ d(y,y_i) \leq \delta_i l_i $, $ R_i >l_i$ and $ R_i-d(y, y_i) \geq R_i-\delta_i l_i>l_i-\delta_i l_i=l_i(1-\delta_i) $, we obtain,

$$ 0 < v_i(z)=\dfrac{u_i(y)}{u_i(y_i)} \leq \dfrac{l_i}{l_i(1-\delta_i)} \leq 2. $$

We set, $ \beta_i= \dfrac{1}{1-\delta_i} $, clearly $ \beta_i \to 1 $.

\bigskip

The function $ v_i $ satisfies the following equation:

$$ -g^{jk}[\exp_{y_i}(y)]\partial_{jk} v_i-\partial_k \left [ g^{jk}\sqrt { |g| } \right ][\exp_{y_i}(y)]\partial_j v_i+ \dfrac{ R_g[\exp_{y_i}(y)]}{[u_i(y_i)]^2} v_i= \tilde V_i{v_i}^3, $$

with, $ \tilde V_i(y)=V_i[\exp_{y_i}(y/[u_i(y_i)])] $. Without loss of generality, we can suppose $ V(x_0)= 8 $.

\bigskip

We use Ascoli and Ladyzenskaya theorems to obtain the uniform convergence (on each compact set of $ {\mathbb R}^4 $) of $ ( v_i)_i $ to $ v $ solution on $ {\mathbb R}^4 $ of: 

$$ \Delta v=8v^3, \,\, v(0)=1,\,\, 0 \leq v\leq 1\leq 2, $$

By the maximum principle, we have $ v>0 $ on $ {\mathbb R}^n $. I we use the Caffarelli-Gidas-Spruck result (see [C-G-S]), we have,$ v(y)= \dfrac{1}{1+{|y|^2}} $.

\bigskip

\underbar {\it Polar Geodesic Coordinates}

\bigskip

Let $ u $ be a function on $ M $. We set $ \bar u(r,\theta)=u[\exp_x(r\theta)] $. We denote $ g_{x,ij} $ the local expression of the metric $ g $ in the exponential chart centered in $ x $.

\bigskip

We set,

$$ w_i(t,\theta)=e^t\bar u_i(e^t,\theta) = e^{t}u_i[\exp_{y_i}(e^t\theta)] \,\,\, {\rm and} \,\,\, \bar V_i(t,\theta)=V_i[\exp_{y_i}(e^t\theta)]. $$

$$ a(y_i,t,\theta)=\log J(y_i,e^t,\theta)=\log [\sqrt { det(g_{y_i, ij})}]. $$ 

We can write the Laplacian in the geodesic polar coordinates:

$$ -\Delta u =\partial_{rr} \bar u+\dfrac{3}{r} \partial_r \bar u+\partial_r [\log  J(x,r,\theta)] \partial_r \bar u-\dfrac{1}{r^2}\Delta_{\theta } \bar u . $$

We deduce the two following lemmas:

\bigskip

\underbar {\bf Lemma 1:}

\smallskip

The function $ w_i $ is a solution of:

$$  -\partial_{tt} w_i-\partial_t a \partial_t w_i-\Delta_{\theta }w_i+c w_i=V_i w_i^3, $$

avec,

 $$ c = c(y_i,t,\theta)= 1+\partial_t a + R_g e^{2t}, $$ 

\underbar {\bf Proof of the Lemma 1:}

\smallskip

We write:

$$ \partial_t w_i=e^{2t}\partial_r \bar u_i+ w_i,\,\, \partial_{tt} w_i=e^{3t} \left [\partial_{rr} \bar u_i+\dfrac{3}{e^t}\partial_r \bar u_i \right ]+ w_i. $$

$$ \partial_t a =e^t\partial_r \log J(y_i,e^t,\theta), \partial_t a \partial_t w_i=e^{3t}\left [ \partial_r \log J\partial_r \bar u_i \right ]+ \partial_t a w_i.$$

Le lemma 1 follows.

\bigskip

Let $ b_1(y_i,t,\theta)=J(y_i,e^t,\theta)>0 $. We can write:

$$ -\dfrac{1}{\sqrt {b_1}}\partial_{tt} (\sqrt { b_1} w_i)-\Delta_{ \theta }w_i+[c(t)+ b_1^{-1/2} b_2(t,\theta)]w_i=\bar V_i {w_i}^3, $$

where, $ b_2(t,\theta)=\partial_{tt} (\sqrt {b_1})=\dfrac{1}{2 \sqrt { b_1}}\partial_{tt}b_1-\dfrac{1}{4(b_1)^{3/2}}(\partial_t b_1)^2 .$

\bigskip

We set,

$$ \tilde w_i=\sqrt {b_1} w_i. $$

\underbar {\bf Lemma 2:}

\smallskip

The function $ \tilde w_i $ is solution of:

$$ -\partial_{tt} \tilde w_i+\Delta_{ \theta } (\tilde w_i)+2\nabla_{\theta}(\tilde  w_i) .\nabla_{\theta} \log (\sqrt {b_1})+(c+b_1^{-1/2} b_2-c_2) \tilde w_i= $$

$$ = \bar V_i \left (\dfrac{1}{b_1} \right )^{1/2} {\tilde w_i}^3, $$

where, $ c_2 $ is a function to be deterined.

\bigskip

\underbar {\bf Proof of the Lemma 2:}

We have: 

$$ -\partial_{tt} \tilde w_i-\sqrt {b_1} \Delta_{ \theta } w_i+(c+b_2) \tilde w_i= \bar V_i \left (\dfrac{1}{b_1} \right )^{1/2} {\tilde w_i}^3, $$

But,

$$ \Delta_{ \theta } (\sqrt {b_1} w_i)=\sqrt {b_1} \Delta_{ \theta } w_i-2 \nabla_{\theta} w_i .\nabla_{\theta} \sqrt {b_1}+ w_i \Delta_{ \theta }(\sqrt {b_1}), $$

and,

$$ \nabla_{\theta} (\sqrt {b_1} w_i)=w_i \nabla_{\theta} \sqrt {b_1}+ \sqrt {b_1} \nabla_{\theta} w_i, $$

we can write,

$$ \nabla_{\theta} w_i. \nabla_{\theta} \sqrt {b_1}=\nabla_{\theta}(\tilde  w_i) .\nabla_{\theta} \log (\sqrt {b_1})-\tilde w_i|\nabla_{\theta} \log (\sqrt {b_1})|^2, $$

we deduce,

$$  \sqrt {b_1} \Delta_{\theta } w_i= \Delta_{ \theta } (\tilde w_i)+2\nabla_{\theta}(\tilde  w_i) .\nabla_{\theta} \log (\sqrt {b_1})-c_2 \tilde w_i, $$

with $ c_2=[\dfrac{1}{\sqrt {b_1}} \Delta_{ \theta }(\sqrt{b_1}) + |\nabla_{\theta} \log (\sqrt {b_1})|^2] . $ The lemma 2 is proved.

\bigskip

\underbar {\it The moving-plane method:}

\bigskip

Let $ \xi_i $  be a real number,  we assume $ \xi_i \leq t $. We set $ t^{\xi_i}=2\xi_i-t $ and $ \tilde w_i^{\xi_i}(t,\theta)=\tilde w_i(t^{\xi_i},\theta) $.

\bigskip

\underbar {\bf Proposition 2:}

\smallskip

We have:

$$ 1)\,\,\, \tilde w_i(\lambda_i,\theta)-\tilde w_i(\lambda_i+4,\theta) \geq \tilde k>0, \,\, \forall \,\, \theta \in {\mathbb S}_{3}. $$

 For all $ \beta >0 $, there exists $ c_{\beta} >0 $ such that:

$$ 2) \,\,\, \dfrac{1}{c_{\beta}} e^t \leq \tilde w_i(\lambda_i+t,\theta) \leq c_{\beta}e^t, \,\, \forall \,\, t\leq \beta, \,\, \forall \,\, \theta \in {\mathbb S}_{3}. $$

\underbar {\bf Proof of the Proposition 2:}

\smallskip

Like in [B], we have, $ w_i(\lambda_i, \theta)-w_i(\lambda_i+4,\theta) \geq k>0 $ for $ i $ large, $ \forall \,\, \theta $. We can remark that $ b_1(y_i,\lambda_i,\theta) \to 1 $ and $ b_1(y_i,\lambda_i+4,\theta) \to 1 $ uniformly in $ \theta $, we obtain 1) of the proposition 2. For 2) we use the previous lemma 2, see also [B].

\bigskip

We set:

$$ \bar Z_i=-\partial_{tt} (...)+\Delta_{ \theta } (...)+2\nabla_{\theta}(...) .\nabla_{\theta} \log (\sqrt {b_1})+(c+b_1^{-1/2} b_2-c_2)(...) $$

{\bf Remark :} In the operator $ \bar Z_i $, we can remark that:

$$ c+b_1^{-1/2}b_2-c_2 \geq k'>0,\,\, {\rm for }\,\, t<<0, $$

it is fundamental if we want to apply the Hopf maximum principle.

\bigskip

\underbar {\bf Goal:}

\bigskip

Like in [B], we have elliptic second order operator. Here it is $ \bar Z_i $, the goal is to use the "moving-plane" method to have a contradiction. For this, we must have:

$$ \bar Z_i(\tilde w_i^{\xi_i}-\tilde w_i) \leq 0, \,\, {\rm if} \,\, \tilde w_i^{\xi_i}-\tilde w_i \leq 0. $$

We write, $ \Delta_{\theta}=\Delta_{g_{y_i, e^t, {}_{{\mathbb S}_{n-1}}}} $. We obtain:

$$ \bar Z_i(\tilde w_i^{\xi_i}-\tilde w_i)= (\Delta_{g_{y_i, e^{t^{\xi_i}}, {}_{{\mathbb S}_{3}}}}-\Delta_{g_{y_i, e^{t}, {}_{{\mathbb S}_{3}}}}) (\tilde w_i^{\xi_i})+ $$

$$ +2(\nabla_{\theta, e^{t^{\xi_i}}}-\nabla_{\theta, e^t})(w_i^{\xi_i}) .\nabla_{\theta, e^{t^{\xi_i}}} \log (\sqrt {b_1^{\xi_i}})+ 2\nabla_{\theta,e^t}(\tilde w_i^{\xi_i}).\nabla_{\theta, e^{t^{\xi_i}}}[\log (\sqrt {b_1^{\xi_i}})-\log \sqrt {b_1}]+ $$ 

$$ +2\nabla_{\theta,e^t} w_i^{\xi_i}.(\nabla_{\theta,e^{t^{\xi_i}}}-\nabla_{\theta,e^t})\log \sqrt {b_1}- [(c+b_1^{-1/2} b_2-c_2)^{\xi_i}-(c+b_1^{-1/2}b_2-c_2)]\tilde w_i^{\xi_i} + $$

$$ + \bar V_i^{\xi_i} \left ( \dfrac{1}{b_1^{\xi_i}} \right )^{1/2} ({\tilde w_i}^{\xi_i})^3-\bar V_i\left ( \dfrac{1}{b_1} \right )^{1/2} {\tilde w_i}^3.\qquad (***1) $$

Clearly, we have:

\bigskip

\underbar {\bf Lemma 3 :}

$$ b_1(y_i,t,\theta)=1-\dfrac{1}{3} Ricci_{y_i}(\theta,\theta)e^{2t}+\ldots, $$

$$ R_g(e^t\theta)=R_g(y_i) + <\nabla R_g(y_i)|\theta > e^t+\dots . $$

According to proposition 1 and lemma 3,

\bigskip

\underbar {\bf Propostion 3 :}

$$ \bar Z_i(\tilde w_i^{\xi_i}-\tilde w_i) \leq |\bar V_i^{\xi_i}-\bar V_i|(b_1^{\xi_i})^{-1/2}(w_i^{\xi_i})^3 + \bar V_i {(b_1^{\xi_i})}^{-1/2}[(\tilde w_i^{\xi_i})^3- \tilde w_i^3] + $$

$$ +C|e^{2t}-e^{2t^{\xi_i}}|\left [|\nabla_{\theta} {\tilde w_i}^{\xi_i}| + |\nabla_{\theta}^2(\tilde w_i^{\xi_i})|+ |Ricci_{y_i}|[\tilde w_i^{\xi_i}+(\tilde w_i^{\xi_i})^3] + | R_g(y_i) | \tilde w_i^{\xi_i} \right ] + C' w_i^{\xi_i} |e^{3t^{\xi_i}}-e^{3t}|. $$

\underbar {\bf Proof of the proposition 3:}

\bigskip

In polar geodesic coordinates (and the Gauss lemma):

$$ g = dt^2+r^2{\tilde g}_{ij}^kd\theta^id\theta^j \,\, {\rm et } \,\, \sqrt { |{\tilde g}^k|}=\alpha^k(\theta) \sqrt {[det(g_{x,ij})]}, $$

where $ \alpha^k $ is the volume element of the unit sphere for the open set $ U^k $.

\smallskip

We can write (with the lemma 3):

$$ |\partial_t b_1(t)|+|\partial_{tt} b_1(t)|+|\partial_{tt} a(t)|\leq C e^{2t}, $$

and,

$$ |\partial_{\theta_j} b_1|+|\partial_{\theta_j,\theta_k} b_1|+\partial_{t,\theta_j}b_1|+|\partial_{t,\theta_j,\theta_k} b_1|\leq C e^{2t}, $$

But,

$$ \Delta_{\theta}=\Delta_{g_{y_i, e^t, {}_{{\mathbb S}_{3}}}}=-\dfrac{\partial_{\theta^l}[{\tilde g}^{\theta^l \theta^j}(e^t,\theta)\sqrt { |{\tilde g}^k(e^t,\theta)|}\partial_{\theta^j}]}{\sqrt {|{\tilde g}^k(e^t,\theta)|}} . $$

Then,

$$ A_i:=\left [{ \left [ \dfrac{\partial_{\theta^l}({\tilde g}^{\theta^l \theta^j}\sqrt { |{\tilde g}^k|}\partial_{\theta^j})}{\sqrt {|{\tilde g}^k|}} \right ] }^{\xi_i}-\left [ \dfrac{\partial_{\theta^l}({\tilde g}^{\theta^l \theta^j}\sqrt { |{\tilde g}^k|}\partial_{\theta^j})}{\sqrt {|{\tilde g}^k|}} \right ] \right ](\tilde w_i^{\xi_i}) = B_i+D_i $$

where,

$$ B_i=\left [ {\tilde g}^{\theta^l \theta^j}(e^{t^{\xi_i}}, \theta)-{\tilde g}^{\theta^l \theta^j}(e^t,\theta) \right ] \partial_{\theta^l \theta^j}\tilde w_i^{\xi_i}, $$

and,

$$ D_i=\left [ \dfrac{\partial_{\theta^l}[{\tilde g }^{\theta^l \theta^j}(e^{t^{\xi_i}},\theta)\sqrt {| {\tilde g}^k|}(e^{t^{\xi_i}},\theta)  ]}{ \sqrt {| {\tilde g}^k|}(e^{t^{\xi_i}},\theta )  }            -\dfrac{\partial_{\theta^l} [{\tilde g }^{\theta^l \theta^j}(e^t,\theta)\sqrt {| {\tilde g}^k|}(e^t,\theta) ]}{ \sqrt {| {\tilde g}^k|}(e^t,\theta) } \right ] \partial_{\theta^j} \tilde w_i^{\xi_i}. $$

Clearly, we can choose $ \epsilon_1 >0 $ such that:

$$ |\partial_r{\tilde g}_{ij}^k(x,r,\theta)|+|\partial_r\partial_{\theta^m}{\tilde g}_{ij}^k(x,r, \theta)| \leq C r,\,\, x\in B(x_0,\epsilon_1) \,\, r\in [0,\epsilon_1], \,\,\theta \in U^k.$$

finally,

$$ A_i \leq C_k|e^{2t}-e^{2t^{\xi_i}}|\left [ |\nabla_{\theta} \tilde w_i^{\xi_i}| + |\nabla_{\theta}^2(\tilde w_i^{\xi_i})| \right ], $$

It is easy to see that: 

$$ \dfrac{ |\nabla_{\theta} ( {\tilde w_i}^{\xi_i}) |}{{\tilde w_i}^{\xi_i}} \leq K \,\,\, {\rm and} \,\,\,  \dfrac{ | \nabla^2_{\theta} ( {\tilde w_i}^{\xi_i} )|}{{\tilde w_i}^{\xi_i}}  \leq K'. $$ 

We take, $ C=\max \{ C_i, 1 \leq i\leq q \} $ and we use $ (***1) $. The proposition 3 is proved.

\bigskip

We have,

$$ c(y_i,t,\theta)= 1 + \partial_t a + R_g e^{2t}, \qquad (\alpha_1) $$ 

$$ b_2(t,\theta)=\partial_{tt} (\sqrt {b_1})=\dfrac{1}{2 \sqrt { b_1}}\partial_{tt}b_1-\dfrac{1}{4(b_1)^{3/2}}(\partial_t b_1)^2 ,\qquad (\alpha_2) $$ 

$$ c_2=[\dfrac{1}{\sqrt {b_1}} \Delta_{ \theta }(\sqrt{b_1}) + |\nabla_{\theta} \log (\sqrt {b_1})|^2], \qquad (\alpha_3) $$

We do a conformal change of the metric such that:

$$ Ricci_{x_0}=R_{\tilde g}(x_0)=0, \sqrt { det ({\tilde g}_{x_0,jk})}=1+O(r^s), s\geq 4, $$

it is given by T. Aubin [Au 1], (see also Lee et Parker, [L,P]).

\bigskip

Without loss of generality, we can assume:

$$ g= \tilde g, \,\,\, R_{\tilde g}(y_i) \to 0 \,\,\, {\rm and} \,\,\,  Ricci_{y_i} \to 0. $$

We assume that $ \lambda \leq \lambda_i+2= - \log u_i(y_i)+2 $.

\bigskip

We work on $ [\lambda,t_i] \times {\mathbb S}_3 $ with $ t_i = \log \sqrt { l_i } \to -\infty $, $ l_i $ as in the proposition 1. For $ i $ large $ \log \sqrt {l_i} >> \lambda_i + 2 $.

\bigskip

The functions  $ v_i $ tend to radially symetric function, then, $ \partial_{\theta_j} w_i^{\lambda} \to 0 $ if $ i \to +\infty $ and,

$$ \dfrac{\partial_{\theta_j}w_i^{\lambda }(t,\theta)}{w_i^{\lambda }}=\dfrac{e^{[(\lambda -\lambda_i)+(\xi_i-t)]} e^{[(\lambda -\lambda_i)+(\xi_i-t)]}(\partial_{\theta_j} v_i)(e^{[(\lambda -\lambda_i)+(\lambda -t)]}\theta)}{e^{[(\lambda-\lambda_i)+(\lambda-t)]}v_i[e^{(\lambda-\lambda_i)+(\lambda - t)}\theta]} \leq {\bar C_i}, $$

where $ \bar C_i $ does not depend on $ \lambda $ and tend to 0. We have also,

$$ |\partial_{\theta} w_i^{\lambda }(t,\theta)|+|\partial_{\theta,\theta} w_i^{\lambda }(t,\theta)|\leq {\tilde C_i} w_i^{\lambda }(t,\theta), \,\, {\tilde C_i} \to 0. $$

and,

$$ |\partial_{\theta} \tilde w_i^{\lambda }(t,\theta)|+|\partial_{\theta,\theta}  \tilde w_i^{\lambda }(t,\theta)|\leq {\tilde C_i} \tilde w_i^{\lambda }(t,\theta), \,\, {\tilde C_i} \to 0. $$

$ \tilde C_i $ does not depend on $ \lambda $.

\bigskip

Now, we set:

$$ \bar w_i=\tilde w_i-\dfrac{\tilde m}{2} e^{t} $$

Like in [B], we have,

\bigskip

\underbar {\bf Lemma 4:}

\bigskip

There is $ \nu <0 $ such that for $ \lambda \leq \nu $ :

$$ \bar w_i^{\lambda}(t,\theta)-\bar w_i(t,\theta) \leq 0, \,\, \forall \,\, (t,\theta) \in [\lambda,t_i] \times {\mathbb S}_{3}. $$

Let $ \xi_i $ be the following real number,

$$ \xi_i=\sup \{ \lambda \leq \lambda_i+2, \bar w_i^{\xi_i}(t,\theta)-\bar w_i(t,\theta) \leq 0, \,\, \forall \,\, (t,\theta)\in [\xi_i,t_i]\times {\mathbb S}_{3} \}. $$

Like in [B], we use the previous lemma to show:

$$ \bar w_i^{\xi_i}-\bar w_i \leq 0 \Rightarrow \bar Z_i(\bar w_i^{\xi_i}-\bar w_i) \leq 0. $$

If we use $ (\alpha_1), (\alpha_2) $ and $ (\alpha_3) $, we have,

$$ \bar Z_i(\bar w_i^{\xi_i}-\bar w_i) \leq 2 A_i (e^t-e^{t^{\xi_i}}) ( \tilde w_i^{\xi_i})^3 + V_i (b_1^{\xi_i})^{-1/2}[( \tilde w_i^{\xi_i})^3- \tilde w_i^3]+ o(1) e^{2t} (e^t-e^{t^{\xi_i}})+o(1) {\tilde w_i}^{\xi_i}(e^{2t}-e^{2t^{\xi_i}}). $$

We can write,

$$ e^{2t}-e^{2t^{\xi_i}}=(e^t-e^{t^{\xi_i}})(e^t+e^{t^{\xi_i}}) \leq 2 e^t (e^t-e^{t^{\xi_i}}). $$

Thus,

$$ \bar Z_i(\bar w_i^{\xi_i}-\bar w_i) \leq 4 e A_i (e^t-e^{t^{\xi_i}}) ({\tilde w_i}^{\xi_i})^2 + V_i (b_1^{\xi_i})^{-1/2}[( \tilde w_i^{\xi_i})^3-{\tilde w_i}^3]+ o(1) e^{2t} (e^t-e^{t^{\xi_i}})+o(1)e^t {\tilde w_i}^{\xi_i}(e^t-e^{t^{\xi_i}}). $$

But,

$$  0 < {\tilde w_i}^{\xi_i} \leq 2e, \,\,\, \tilde w_i \geq \dfrac{m}{2} e^t \,\,\, {\rm and} \,\,\, {\tilde w_i}^{\xi_i}-{\tilde w_i} \leq \dfrac{m}{2} ( e^{t^{\xi_i}}-e^t), $$

and,

$$ ({\tilde w_i}^{\xi_i})^3-{\tilde w_i}^3=({\tilde w_i}^{\xi_i}-\tilde w_i)[ ({\tilde w_i}^{\xi_i})^2+{\tilde w_i}^{\xi_i} {\tilde w_i}+ {\tilde w_i}^2] \leq ({\tilde w_i}^{\xi_i}-{\tilde w_i}) ({\tilde w_i^{\xi_i}})^2 + ({\tilde w_i}^{\xi_i}-{\tilde w_i}) \dfrac{m^2 e^{2t}}{4} + ({\tilde w_i}^{\xi_i}-{\tilde w_i}) \dfrac{m}{2} e^t {\tilde w_i}^{\xi_i}, $$

then,

$$ \bar Z_i(\bar w_i^{\xi_i}-\bar w_i) \leq    \left [ ({\tilde w_i}^{\xi_i})^2 [\dfrac{a m}{4} - 4 e A_i] + [\dfrac{am^3}{16}-o(1)] + [\dfrac{am^2}{8}- o(1)] e^t {\tilde w_i}^{\xi_i} \right ] ( e^{t^{\xi_i}}- e^t) \leq 0. $$

I fwe use the Hopf maximum principle, we obtain (like in [B]):

$$ \max_{\theta \in {\mathbb S}_3} w_i(t_i,\theta) \leq \min_{\theta \in {\mathbb S}_3} w_i(2\xi_i-t_i), $$

we can write (by using the proposition 2):

$$  l_i u_i(y_i) \leq c, $$

Contradiction.

\bigskip

\underbar{\bf References:}

{\small \bigskip

[Au 1] T. Aubin. Some Nonlinear Problems in Riemannian Geometry. Springer-Verlag 1998.

\smallskip

[Au 2] T. Aubin. Sur quelques probl\`emes de courbure scalaire in J. Func. Anal 2006.

\smallskip

[B] S.S Bahoura. Majorations du type $ \sup u \times \inf u \leq c $ pour l'\'equation de la courbure scalaire sur un ouvert de $ {\mathbb R}^n, n\geq 3 $. J. Math. Pures. Appl.(9) 83 2004 no, 9, 1109-1150.

\smallskip

[B-L-S] H. Brezis, YY. Li, I. Shafrir. A sup+inf inequality for some
nonlinear elliptic equations involving exponential
nonlinearities. J.Funct.Anal.115 (1993) 344-358.

\smallskip

[C-G-S] L. Caffarelli, B. Gidas, J. Spruck. Asymptotic symmetry and local
behavior of semilinear elliptic equations with critical Sobolev
growth. Comm. Pure Appl. Math. 37 (1984) 369-402.

\smallskip

[C-L] C-C.Chen, C-S. Lin. Estimates of the conformal scalar curvature
equation via the method of moving planes. Comm. Pure
Appl. Math. L(1997) 0971-1017.

\smallskip

[D] O. Druet. Compactness for Yamabe metrics in low diemensions, Int. Math. Res. Not. 23 (2004) 1143-1191.

\smallskip

[L,P] J.M. Lee, T.H. Parker. The Yamabe problem. Bull.Amer.Math.Soc (N.S) 17 (1987), no.1, 37 -91.

\smallskip

[L] YY. Li. Prescribing scalar curvature on $ {\mathbb S}_n $ and related
Problems. C.R. Acad. Sci. Paris 317 (1993) 159-164. Part
I: J. Differ. Equations 120 (1995) 319-410. Part II: Existence and
compactness. Comm. Pure Appl.Math.49 (1996) 541-597.

\smallskip

[L-Z 1] YY. Li, L. Zhang. A Harnack type inequality for the Yamabe equation in low dimensions.  Calc. Var. Partial Differential Equations  20  (2004),  no. 2, 133--151

\smallskip

[L-Z 2] YY. Li, L. Zhang. Compactness of solutions to the Yamabe problem. II.  Calc. Var. Partial Differential Equations  24  (2005),  no. 2, 185--237. 

\smallskip

[L-Zh] YY. Li, M. Zhu. Yamabe type equations on three-dimensional Riemannian manifolds.  Commun. Contemp. Math.  1  (1999),  no. 1, 1--50. 

\smallskip

[M] F.C. Marques. A priori estimates for the Yamabe problem in the non-locally conformally flat case. J. Differential Geom.  71  (2005),  no. 2, 315--346.

\end{document}